\documentclass[11pt]{article}
\usepackage{amscd}
\usepackage{amsfonts}
\usepackage{amsmath}
\usepackage{amssymb}
\usepackage{amsthm}
\usepackage{bbm}
\usepackage{CJK}
\usepackage{fancyhdr}
\usepackage{graphicx}
\usepackage{hyperref}
\usepackage{indentfirst}
\usepackage{latexsym}
\usepackage{mathrsfs}
\usepackage{xypic}
\usepackage{amsmath,amssymb,amsfonts,amsthm,fancyhdr}
\usepackage{epsfig,graphicx,picins,picinpar,subfigure}
\usepackage{pstricks}
\usepackage{amssymb}
\usepackage{xypic}
\usepackage[compress]{cite}

\newtheorem{theorem}{Theorem}[section]
\newtheorem{prop}[theorem]{Proposition}
\newtheorem{lemma}[theorem]{Lemma}
\newtheorem{coro}[theorem]{Corollary}
\newtheorem{prop-def}{Proposition-Definition}[section]

\theoremstyle{definition}
\newtheorem{defn}[theorem]{Definition}

\newtheorem{remark}[theorem]{Remark}
\newtheorem{exam}[theorem]{Example}

\usepackage[top=1in,bottom=1in,left=1.25in,right=1.25in]{geometry}
\def\<{\langle}
\def\>{\rangle}

\date{August 1, 2023}
\begin{document}
\renewcommand{\baselinestretch}{1.2}
\renewcommand{\arraystretch}{1.0}
\title{\bf   Deformations and   extensions of  modified $\lambda$-differential $3$-Lie Algebras}
\author{{\bf Wen Teng$^{1}$,   Hui Zhang$^{2,3,*}$}\\
{\small 1. School of Mathematics and Statistics, Guizhou University of Finance and Economics,} \\
{\small  Guiyang  550025, P. R. of China}\\
  {\small Email: tengwen@mail.gufe.edu.cn (Wen Teng)} \\
{\small 2.   School of Information, Guizhou University of Finance and Economics,   }\\
{\small Guiyang  $550005$, P. R. of China}\\
{\small 3. Postdoctoral Scientific Research Station, ShijiHengtong Technology Co.,} \\
{\small   Ltd, Guiyang,  550014, P. R. of China}\\
  {\small * Corresponding author: zhanghui@mail.gufe.edu.cn (Hui Zhang)} \\}
 \maketitle
\begin{center}
\begin{minipage}{13.cm}

{\bf Abstract}
In this paper, we introduce the  representation of modified $\lambda$-differential $3$-Lie algebras and   define the cohomology of modified $\lambda$-differential $3$-Lie algebras with coefficients in a   representation.  As   applications of the proposed cohomology theory,  we study  linear  deformations, abelian extensions and $T^*$-extensions of modified $\lambda$-differential $3$-Lie algebras.
 \smallskip

 \smallskip

 \end{minipage}
 \end{center}
 \normalsize\vskip0.5cm

\section{Introduction}
\def\theequation{\arabic{section}. \arabic{equation}}
\setcounter{equation} {0}

3-Lie algebra plays an important role in string theory, and it is also used to study supersymmetry and gauge symmetry transformation of the world-volume theory of  multiple coincident M2-branes \cite{Bagger,Gustavsson}.
The concept of 3-Lie algebra,    general $n$-Lie algebra, was first introduced by Filippov \cite{Filippov}, which can be regarded as a  generalization of Lie algebra to higher-order algebra.
 3-Lie algebra has attracted the attention of scholars from mathematics and physics \cite{Cherkis,Ho, Papadopoulos}. Representation theory, cohomology theory, deformations, Nijenhuis operators and extension theory  of $n$-Lie algebras have been widely studied by scholars \cite{Kasymov,Takhtajan,Arfa,Mignel, Liu, Liu1,LiuW,Sheng, Xu, Zhang,Zhao}.

  Derivations  play an important roles in the study of homotopy Lie algebras \cite{Voronov},  differential Galois theory \cite{Magid}, control theory and gauge theories of quantum field theory \cite{Ayala1}.
Recently,   authors  have studied  algebras with derivations in \cite{Doubek,Loday} from the
operadic point of view.
 In \cite{Tang},   authors have studied the deformation and extension of Lie algebras with derivations  from the cohomological point of view.
The results of \cite{Tang}   have been extended to 3-Lie algebras with derivations \cite{GuoSaha,XuLiu}. More research on algebraic structures with derivations  have
been developed, see \cite{Das1,Wu1,Wu2,Sun, SunWu, Guo4} and references cited therein.

In recent years, due to the important work of \cite{Bai,Guo0,Guo1,Guo2,Das2, Wang},   scholars   pay attention to the structures with arbitrary weights.
For $\lambda\in  \mathbf{k}$,  the cohomology, extensions and
deformations of  Rota-Baxter 3-Lie algebras with any weight $\lambda$ and differential 3-Lie algebras with any weight $\lambda$  were established in  \cite{Hou,Guo3,SunChen}. In addition, in \cite{Das4},   author has  studied the cohomology and deformation   of modified Rota-Baxter  algebras.  The cohomology and deformation   of modified Rota-Baxter Leibniz algebras with weight $\lambda$ are given in \cite{Li1,Mondal}.
In \cite{Peng},   the concept  of  modified   $\lambda$-differential Lie algebras has been introduced.
Motivated by the mentioned
work on the modified $\lambda$-differential operator
on Lie algebras,
in this paper, our main objective is  to  study  modified $\lambda$-differential 3-Lie algebras. We develop
a cohomology theory of modified $\lambda$-differential 3-Lie algebras    that controls the deformations  and extensions   of  modified $\lambda$-differential 3-Lie algebras.

The paper is organized as follows.   In Section  \ref{sec:Repre}, we  introduce the  representation  and  cohomology  of    modified $\lambda$-differential  $3$-Lie algebras.
 In Section \ref{sec:def},  we   consider linear deformations  of    modified $\lambda$-differential 3-Lie algebras.
 In Section \ref{sec:aext},  we study abelian extensions   of   modified $\lambda$-differential 3-Lie algebras.
 In Section \ref{sec:Text},  we study $T^*$-extensions   of   modified $\lambda$-differential 3-Lie algebras.

Throughout this paper, $ \mathbf{k}$ denotes a field of characteristic zero. All the   vector spaces and
  linear maps  are taken over $ \mathbf{k}$.

\section{Representations and cohomologies  of modified $\lambda$-differential $3$-Lie algebras}\label{sec:Repre}
\def\theequation{\arabic{section}.\arabic{equation}}
\setcounter{equation} {0}

In this section, first, we  introduce the  concept of  a modified $\lambda$-differential $3$-Lie algebras.
Then, we  give representations and cohomologies of    modified $\lambda$-differential  $3$-Lie algebras.

\begin{defn}  \cite{Filippov}
(i) A $3$-Lie algebra is a    tuple $(\mathfrak{A}, [-, -, -])$ in which $\mathfrak{A}$ is a vector space together with  a skew-symmetric ternary operation $[-, -, -]:\wedge^3\mathfrak{A}\rightarrow \mathfrak{A}$   such that
\begin{align}
 &[a_1, a_2, [a_3, a_4, a_5]]=[[a_1, a_2, a_3], a_4, a_5]+ [a_3,  [a_1, a_2,  a_4], a_5]+ [a_3, a_3,  [a_1, a_2,  a_5]],\label{2.1}
\end{align}
for all $ a_1, a_2,  a_3, a_4, a_5\in \mathfrak{A}$.\\
(ii) A homomorphism between two  $3$-Lie algebras  $(\mathfrak{A}_1, [-, -, -]_1)$ and $(\mathfrak{A}_2, [-, -,-]_2)$ is a linear map $\eta: \mathfrak{A}_1\rightarrow \mathfrak{A}_2$ satisfying
$\eta([a_1, a_2,  a_3]_1)=[\eta(a_1), \eta(a_2),\eta(a_3)]_2,~~\forall ~a_1, a_2,  a_3\in \mathfrak{A}_1.$
\end{defn}

\begin{defn}
(i)  Let  $\lambda\in \mathbf{k}$ and  $(\mathfrak{A},  [-, -, -])$ be a  3-Lie algebra. A modified $\lambda$-differential operator on 3-Lie algebra $\mathfrak{A}$   is a linear operator $\mathrm{d}:\mathfrak{A}\rightarrow \mathfrak{A}$,  such that
\begin{align}
\mathrm{d}[a_1, a_2,  a_3]=[\mathrm{d}(a_1), a_2,  a_3]+[a_1,  \mathrm{d}(a_2),  a_3]+[a_1, a_2,  \mathrm{d}(a_3)]+\lambda[a_1, a_2,  a_3],\forall a_1, a_2,  a_3\in \mathfrak{A}.\label{2.2}
\end{align}
(ii) A modified $\lambda$-differential  3-Lie algebra   is a triple  $(\mathfrak{A}, [-, -,-], \mathrm{d})$  consisting of a
3-Lie algebra  $(\mathfrak{A}, [-, -,-])$ and a modified $\lambda$-differential operator $\mathrm{d}$.\\
(iii) A homomorphism between two modified $\lambda$-differential  3-Lie algebras  $(\mathfrak{A}_1, [-, -, -]_1,\mathrm{d}_1)$ and $(\mathfrak{A}_2,   [-,-, -]_2, \mathrm{d}_2)$  is a 3-Lie algebra  homomorphism $\eta: (\mathfrak{A}_1,  [-, -, -]_1)\rightarrow (\mathfrak{A}_2,  [-,-, -]_2)$ such that
$\eta\circ \mathrm{d}_1=\mathrm{d}_2\circ\eta$.  Furthermore, if $\eta$ is  nondegenerate,  then $\eta$ is called an isomorphism from $\mathfrak{A}_1$ to $\mathfrak{A}_2$.
\end{defn}

Let $(\mathfrak{A}, [-,-,-])$ be a 3-Lie algebra, then the elements in $\wedge^2 \mathfrak{A}$ are called fundamental objects of the 3-Lie algebra $(\mathfrak{A}, [-,-,-])$. There is a bilinear operation $[-,-]$ on $\wedge^2 \mathfrak{A}$, which is given by
\begin{eqnarray*}
&&[A,B]_\mathcal{F}=[a_1,a_2,b_1]\wedge  b_2 + b_1 \wedge[a_1,a_2,b_2], \forall A=a_1\wedge a_2, B=b_1\wedge b_2\in  \wedge^2  \mathfrak{A}.
\end{eqnarray*}
 It is shown in   \cite{Daletskii} that $(\wedge^2 \mathfrak{A}, [-,-]_\mathcal{F})$ is a Leibniz algebra. Furthermore, we have the following proposition.

\begin{prop}
Let $(\mathfrak{A}, [-, -, -],\mathrm{d})$  be a  modified $\lambda$-differential 3-Lie algebra. Then, $(\wedge^2 \mathfrak{A}, [-,-]_\mathcal{F}, \mathrm{d}_{\mathcal{F}})$ is a    Leibniz algebra with a derivation, where $\mathrm{d}_{\mathcal{F}}(a_1\wedge a_2)=\mathrm{d}(a_1)\wedge a_2+a_1\wedge\mathrm{d}(a_2)+\lambda a_1\wedge a_2$, for all $a_1\wedge a_2\in  \wedge^2  \mathfrak{A}$. See \cite{Das1} for more
details about Leibniz algebras with derivations
\end{prop}

\begin{proof}
For any $A=a_1\wedge a_2, B=b_1\wedge b_2\in  \wedge^2  \mathfrak{A}$, by Eq.  \eqref{2.2}  we have
\begin{align*}
&\mathrm{d}_\mathcal{F}[A,B]_\mathcal{F}\\
=&\mathrm{d}[a_1,a_2,b_1]\wedge  b_2 +[a_1,a_2,b_1]\wedge  \mathrm{d}(b_2)+\lambda [a_1,a_2,b_1]\wedge  b_2 \\
&+ \mathrm{d}(b_1) \wedge[a_1,a_2,b_2]+b_1 \wedge\mathrm{d}[a_1,a_2,b_2]+\lambda b_1 \wedge[a_1,a_2,b_2]\\
=&([\mathrm{d}(a_1), a_2,  b_1]+[a_1,  \mathrm{d}(a_2), b_1]+[a_1, a_2,  \mathrm{d}(b_1)]+\lambda[a_1, a_2,  b_1])\wedge  b_2 \\
&+[a_1,a_2,b_1]\wedge  \mathrm{d}(b_2)+\lambda [a_1,a_2,b_1]\wedge  b_2 + \mathrm{d}(b_1) \wedge[a_1,a_2,b_2]\\
&+b_1 \wedge([\mathrm{d}(a_1), a_2,  b_2]+[a_1,  \mathrm{d}(a_2), b_2]+[a_1, a_2,  \mathrm{d}(b_2)]+\lambda[a_1, a_2,  b_2])+\lambda b_1 \wedge[a_1,a_2,b_2]\\
=&([\mathrm{d}(a_1), a_2,  b_1]+[a_1,  \mathrm{d}(a_2), b_1]+\lambda[a_1, a_2,  b_1])\wedge  b_2+b_1 \wedge([\mathrm{d}(a_1), a_2,  b_2]+[a_1,  \mathrm{d}(a_2), b_2]+\lambda[a_1, a_2,  b_2])\\
&+[a_1, a_2,  \mathrm{d}(b_1)])\wedge  b_2+[a_1,a_2,b_1]\wedge  \mathrm{d}(b_2) +\lambda [a_1,a_2,b_1]\wedge  b_2\\
&+b_1\wedge[a_1, a_2,  \mathrm{d}(b_2)]+ \mathrm{d}(b_1) \wedge[a_1,a_2,b_2]+\lambda b_1 \wedge[a_1,a_2,b_2]\\
=&[\mathrm{d}_\mathcal{F}(A),B]_\mathcal{F}+[A,\mathrm{d}_\mathcal{F}(B)]_\mathcal{F}.
\end{align*}
Hence, $(\wedge^2 \mathfrak{A}, [-,-]_\mathcal{F}, \mathrm{d}_{\mathcal{F}})$ is a    Leibniz algebra with a derivation.
\end{proof}

\begin{remark}
Let $\mathrm{d}$ be a modified $\lambda$-differential operator on  3-Lie algebra  $(\mathfrak{A},  [-, -, -])$. If $\lambda=0,$ then $\mathrm{d}$ is a derivation on 3-Lie algebra  $\mathfrak{A}$.  See \cite{ZhangK} for   various derivations of 3-Lie algebras.
\end{remark}

\begin{exam}
Let   $(\mathfrak{A}, [-, -, -])$ be a  3-Lie algebra.  Then, a linear operator $\mathrm{d}:\mathfrak{A}\rightarrow \mathfrak{A}$ is a   modified $\lambda$-differential operator
if and only if $\mathrm{d}+\frac{\lambda}{2} \mathrm{id}_\mathfrak{A}$ is a derivation on 3-Lie algebra $\mathfrak{A}$.
\end{exam}

\begin{exam}
Let $(\mathfrak{A},[-,-, -],\mathrm{d})$ be a modified $\lambda$-differential   3-Lie algebra.
Then, for $k\in \mathbf{k}$, $(\mathfrak{A},[-,-,-],k\mathrm{d})$ is a modified $(k\lambda)$-differential   3-Lie algebra.
\end{exam}

\begin{exam}
Let $(\mathfrak{A},  [-, -, -])$ be a 3-dimensional   3-Lie algebra   with a basis $\mathfrak{e}_1$, $\mathfrak{e}_2$ and $\mathfrak{e}_3$ defined by
$$ [\mathfrak{e}_1,\mathfrak{e}_2,\mathfrak{e}_3]=\mathfrak{e}_1.$$
Then, any modified  $(-k_{22}-k_{33})$-differential  operator $\mathrm{d}=(k_{ij})$ has the form as follows
$$\mathrm{d}(\mathfrak{e}_1,\mathfrak{e}_2,\mathfrak{e}_3)=(\mathfrak{e}_1,\mathfrak{e}_2,\mathfrak{e}_3)\left(
        \begin{array}{ccc}
          k_{11} & k_{12} & k_{13}  \\
          0 & k_{22}  & k_{23}  \\
          0 & k_{32} & k_{33}\\
        \end{array}
      \right),$$
for $k_{ij}\in  \mathbf{k},(i,j=1,2,3)$.
\end{exam}


\begin{defn} (i) (see \cite{Kasymov})
A representation of  a  3-Lie algebra $(\mathfrak{A}, [-, -,-])$ on  a vector space $\mathfrak{M}$ is a skew-symmetric bilinear map $\rho: \mathfrak{A}\wedge \mathfrak{A}\rightarrow \mathrm{End}(\mathfrak{M})$, such that
\begin{align}
& \rho([a_1, a_2,  a_3],a_4)=\rho(a_2,  a_3)\rho(a_1,a_4)+\rho(a_3,  a_1)\rho(a_2,  a_4)+\rho(a_1,  a_2)\rho(a_3,  a_4),\label{2.3}\\
& \rho(a_1,  a_2)\rho(a_3,  a_4)=\rho(a_3,  a_4)\rho(a_1,  a_2)+ \rho([a_1, a_2,  a_3],a_4)+ \rho(a_3,[a_1, a_2,  a_4]),\label{2.4}\
\end{align}
for all  $a_1, a_2,  a_3,a_4\in \mathfrak{A}$.  We also denote a representation of $ \mathfrak{A} $  on   $\mathfrak{M}$  by  $(\mathfrak{M};  \rho)$. \\
(ii) A representation of a modified $\lambda$-differential 3-Lie algebra  $(\mathfrak{A}, [-, -,-],\mathrm{d})$   is a triple $(\mathfrak{M};  \rho, \mathrm{d}_\mathfrak{M})$, where $(\mathfrak{M};  \rho)$ is a representation of the 3-Lie algebra  $(\mathfrak{A}, [-, -, -])$  and  $\mathrm{d}_\mathfrak{M}$  is a linear operator on $\mathfrak{M}$, satisfying the following equation
\begin{align}
\mathrm{d}_\mathfrak{M}(\rho(a,b)v)=\rho(\mathrm{d}(a),b)v+\rho(a,\mathrm{d}(b))v+\rho(a,b)\mathrm{d}_\mathfrak{M}(v)+\lambda\rho(a,b)v, \label{2.5}
\end{align}
for any $a,b\in \mathfrak{A}$ and $ v\in \mathfrak{M}.$
\end{defn}

\begin{remark}
Let $(\mathfrak{M};  \rho, \mathrm{d}_\mathfrak{M})$ be a representation of the modified $\lambda$-differential 3-Lie algebra $(\mathfrak{A}, [-, -, -],\mathrm{d})$. If $\lambda=0,$ then $(\mathfrak{M};  \rho, \mathrm{d}_\mathfrak{M})$ is a representation of the 3-Lie algebra with a derivation $(\mathfrak{A}, [-, -, -],\mathrm{d})$.
One can refer to \cite{GuoSaha,XuLiu,SunWu} for more information about 3-Lie algebras with  derivations.
\end{remark}

\begin{exam}
Let $(\mathfrak{M};  \rho)$ be a representation of the  3-Lie algebra $(\mathfrak{L}, [-, -, -])$.  Then $(\mathfrak{M};  \rho, \mathrm{d}_\mathfrak{M})$ is a representation of the modified $\lambda$-differential 3-Lie algebra $(\mathfrak{A}, [-, -, -], \mathrm{d})$ if and only if $(\mathfrak{M};  \rho, \mathrm{d}_\mathfrak{M}+\frac{\lambda}{2} \mathrm{id}_\mathfrak{M})$ is a representation of the  3-Lie algebra with a derivation $(\mathfrak{A}, [-, -, -],\mathrm{d}+\frac{\lambda}{2} \mathrm{id}_\mathfrak{A})$.
\end{exam}

\begin{exam}
Let $(\mathfrak{M};  \rho)$ be a representation of  the 3-Lie algebra $(\mathfrak{A}, [-, -, -])$.
Then, for $k\in  \mathbf{k}$, $(\mathfrak{M}; \rho, \mathrm{id}_\mathfrak{M})$ is a representation of the modified $(-2k)$-differential 3-Lie algebra $(\mathfrak{A}, [-, -, -],k\mathrm{id}_\mathfrak{A})$.
\end{exam}

\begin{exam}
Let $(\mathfrak{M};  \rho, \mathrm{d}_\mathfrak{M})$ be a representation of the modified $\lambda$-differential 3-Lie algebra $(\mathfrak{A}, [-, -, -],\mathrm{d})$.
Then, for $k\in \mathbf{k}$, $(\mathfrak{M};  \rho, k \mathrm{d}_\mathfrak{M})$ is a representation of the modified $(k\lambda)$-differential 3-Lie algebra $(\mathfrak{A}, [-, -, -],k \mathrm{d})$.
\end{exam}

\begin{prop}
Let   $(\mathfrak{A}, [-, -, -])$ be a  3-Lie algebra, and  $(\mathfrak{M};  \rho)$ be a   representation of it. Then  $(\mathfrak{M};  \rho, \mathrm{d}_\mathfrak{M})$ is a representation of modified $\lambda$-differential 3-Lie algebra  $(\mathfrak{A}, [-, -, -],\mathrm{d})$ if and only if  $\mathfrak{A }\oplus \mathfrak{M}$ is a  modified $\lambda$-differential 3-Lie algebra under the following maps:
\begin{align*}
[a_1+u_1, a_2+u_2, a_3+u_3]_{\rho}:=&[a_1, a_2, a_3]+\rho(a_1,a_2)u_3+\rho(a_3,a_1)u_2+\rho(a_2,a_3)u_1,\\
\mathrm{d} \oplus \mathrm{d}_\mathfrak{M}(a_1+u_1):=&\mathrm{d}(a_1)+\mathrm{d}_\mathfrak{M}(u_1),
\end{align*}
for all  $a_1, a_2, a_3\in \mathfrak{A}$ and $u_1, u_2, u_3\in \mathfrak{M}$.
\end{prop}

\begin{proof}
At first, it is known that $(\mathfrak{A} \oplus \mathfrak{M},[-, -, -]_{\rho})$    is a 3-Lie algebra.
  Next, for any $a_1,a_2,a_3\in \mathfrak{A}, u_1,u_2,u_3\in \mathfrak{M}$, by Eqs. \eqref{2.2}   and \eqref{2.5}, we have
\begin{align*}
&\mathrm{d} \oplus \mathrm{d}_\mathfrak{M}([a_1+u_1,a_2+u_2,a_3+u_3]_{\rho})\\
=&\mathrm{d} [a_1, a_2, a_3]+ \mathrm{d}_\mathfrak{M}(\rho(a_1,a_2)u_3)+\mathrm{d}_\mathfrak{M}(\rho(a_3,a_1)u_2)+\mathrm{d}_\mathfrak{M}(\rho(a_2,a_3)u_1)\\
=&[\mathrm{d}(a_1), a_2,  a_3]+[a_1,  \mathrm{d}(a_2),  a_3]+[a_1, a_2, \mathrm{d}(a_3)]+\lambda[a_1, a_2,  a_3]\\
&+\rho(\mathrm{d}(a_1),a_2)u_3+\rho(a_1,\mathrm{d}(a_2))u_3+\rho(a_1,a_2)\mathrm{d}_\mathfrak{M}(u_3)+\lambda\rho(a_1,a_2)u_3\\
&+\rho(\mathrm{d}(a_3),a_1)u_2+\rho(a_3,\mathrm{d}(a_1))u_2+\rho(a_3,a_1)\mathrm{d}_{\mathfrak{M}}(u_2)+\lambda\rho(a_3,a_1)u_2\\
&+\rho(\mathrm{d}(a_2),a_3)u_1+\rho(a_2,\mathrm{d}(a_3))u_1+\rho(a_2,a_3)\mathrm{d}_{\mathfrak{M}}(u_1)+\lambda\rho(a_2,a_3)u_1\\
=&[\mathrm{d} \oplus \mathrm{d}_\mathfrak{M}(a_1+u_1),a_2+u_2,a_3+u_3]_{\rho}+[a_1+u_1,\mathrm{d} \oplus \mathrm{d}_\mathfrak{M}(a_2+u_2),a_3+u_3]_{\rho}\\
&+[a_1+u_1,a_2+u_2,\mathrm{d }\oplus \mathrm{d}_\mathfrak{M}(a_3+u_3)]_{\rho}+\lambda[a_1+u_1,a_2+u_2,a_3+u_3]_{\rho}.
\end{align*}
Hence, $(\mathfrak{A} \oplus \mathfrak{M},[-, -, -]_{\rho},\mathrm{d} \oplus \mathrm{d}_\mathfrak{M})$ is a modified $\lambda$-differential 3-Lie algebra.

Conversely, assume that $(\mathfrak{A }\oplus \mathfrak{M},[-,-,-]_{\rho},\mathrm{d} \oplus \mathrm{d}_\mathfrak{M})$ is a  modified $\lambda$-differential 3-Lie algebra,
for any $a_1, a_2\in \mathfrak{A}$ and $  u_3\in \mathfrak{M}$, we have
\begin{align*}
&\mathrm{d} \oplus \mathrm{d}_\mathfrak{M}[a_1+0, a_2+0, 0+u_3]_{\rho}\\
=&[\mathrm{d} \oplus \mathrm{d}_\mathfrak{M}(a_1+0), a_2+0, 0+u_3]_{\rho}+[a_1+0, \mathrm{d} \oplus \mathrm{d}_\mathfrak{M}(a_2+0), 0+u_3]_{\rho}\\
&+[a_1+0, a_2+0, \mathrm{d} \oplus \mathrm{d}_\mathfrak{M}(0+u_3)]_{\rho}+\lambda[a_1+0, a_2+0, 0+u_3]_{\rho},
\end{align*}
which implies that,
\begin{align*}
\mathrm{d}_\mathfrak{M}(\rho(a_1,a_2)u_3)=\rho(\mathrm{d}(a_1),a_2)u_3+\rho(a_1,\mathrm{d}(a_2))u_3+\rho(a_1,a_2)\mathrm{d}_\mathfrak{M}(u_3)+\lambda\rho(a_1,a_2)u_3.
\end{align*}
Therefore, $(\mathfrak{M};  \rho, \mathrm{d}_\mathfrak{M})$ is a representation of    $(\mathfrak{A}, [-, -, -],\mathrm{d})$.
\end{proof}

Let $(\mathfrak{M};  \rho, \mathrm{d}_\mathfrak{M})$  be a representation of a modified $\lambda$-differential  3-Lie algebra $(\mathfrak{A}, [-, -, -],\mathrm{d})$, and $\mathfrak{M}^*:=\mathrm{Hom}(\mathfrak{M},\mathbf{k})$ be a dual space of $\mathfrak{M}$. We define a bilinear map $\rho^*: \wedge^2\mathfrak{A}\rightarrow \mathrm{End}(\mathfrak{M}^*)$ and a linear map $\mathrm{d}_\mathfrak{M}^*: \mathfrak{M}^*\rightarrow \mathfrak{M}^*$, respectively by
\begin{align}
&&\langle \rho^*(a_1,a_2)u^*,v \rangle=-\langle u^*,\rho(a_1,a_2)v \rangle, ~\mathrm{and}~\langle\mathrm{ d}_\mathfrak{M}^*u^*,v \rangle=\langle u^*,\mathrm{d}_\mathfrak{M}(v) \rangle, \label{2.6}
\end{align}
for any $a_1,a_2\in \mathfrak{A}, v\in \mathfrak{M}$ and $u^*\in \mathfrak{M}^*.$

\begin{prop}
 With the above notations, $(\mathfrak{M}^*; \rho^*,-\mathrm{d}^*_\mathfrak{M})$ is a representation of modified $\lambda$-differential 3-Lie algebra $(\mathfrak{A}, [-, -, -],\mathrm{d})$.

\end{prop}
\begin{proof}
First,   It has been proved that \cite{ Rotkiewicz}   $(\mathfrak{M}^*; \rho^*)$ is a representation of the 3-Lie algebra $(\mathfrak{A}, [-, -, -])$.
Furthermore, for any $a_1,a_2\in \mathfrak{A},v\in \mathfrak{M}$ and $u^*\in \mathfrak{M}^*,$ by Eqs.  \eqref{2.5} and \eqref{2.6}, we have
\begin{align*}
&\langle \rho^*(\mathrm{d}(a_1),a_2)u^*,v \rangle+\langle\rho^*(a_1,\mathrm{d}(a_2))u^*,v \rangle+\langle\rho^*(a_1,a_2)(-\mathrm{d}^*_\mathfrak{M})u^*,v \rangle+\langle\lambda\rho^*(a_1,a_2)u^*,v \rangle\\
&- \langle(-\mathrm{d}^*_\mathfrak{M})\rho^*(a_1,a_2)u^*,v \rangle\\
=&-\langle u^*, \rho(\mathrm{d}(a_1),a_2)v \rangle-\langle u^*,\rho(a_1,\mathrm{d}(a_2))v \rangle-\langle (-\mathrm{d}^*_\mathfrak{M})u^*, \rho(a_1,a_2)v  \rangle-\langle u^*,\lambda\rho(a_1,a_2)v \rangle\\
&+\langle \rho^*(a_1,a_2)u^*,\mathrm{d}_\mathfrak{M}(v) \rangle\\
=&-\langle u^*, \rho(\mathrm{d}(a_1),a_2)v \rangle-\langle u^*,\rho(a_1,\mathrm{d}(a_2))v \rangle+\langle u^*,\mathrm{d}_\mathfrak{M}(\rho(a_1,a_2)v) \rangle-\langle u^*,\lambda\rho(a_1,a_2)v \rangle\\
&-\langle u^*,\rho(a_1,a_2)\mathrm{d}_\mathfrak{M}(v) \rangle\\
=&-\langle u^*, \rho(\mathrm{d}(a_1),a_2)v +\rho(a_1,\mathrm{d}(a_2))v-\mathrm{d}_\mathfrak{M}(\rho(a_1,a_2)v) +\lambda\rho(a_1,a_2)v+\rho(a_1,a_2)\mathrm{d}_\mathfrak{M}(v) \rangle\\
=&0,
\end{align*}
which implies that $ \rho^*(\mathrm{d}(a_1),a_2)u^*+\rho^*(a_1,\mathrm{d}(a_2))u^*+\rho^*(a_1,a_2)(-\mathrm{d}^*_\mathfrak{M})u^*+\lambda\rho^*(a_1,a_2)u^*-(-\mathrm{d}^*_\mathfrak{M})\rho^*(a_1,a_2)u^*=0$.
So we get the result.
\end{proof}

\begin{exam} \label{exam:dual rep}
Let  $(\mathfrak{A}, [-, -, -],\mathrm{d})$ be a   modified $\lambda$-differential 3-Lie algebra and define $ ad: \mathfrak{A}\wedge \mathfrak{A} \rightarrow \mathrm{End}(\mathfrak{A})$ by $ ad(a_1,a_2)(a)= [a_1,a_2,a], \forall a_1,a_2,a\in \mathfrak{A}$.
Then    $(\mathfrak{A};  ad,  \mathrm{d})$ is a representation of the modified $\lambda$-differential 3-Lie algebra $(\mathfrak{A}, [-, -, -],  \mathrm{d})$,
 which is called the adjoint representation of $(\mathfrak{A}, [-, -, -],  \mathrm{d})$.
Furthermore,
 $(\mathfrak{A}^*;    ad^* , -\mathrm{d}^*)$ is a dual adjoint representation of  $(\mathfrak{A}, [-, -, -],\mathrm{d})$.
\end{exam}

In the next, we  will study  the cohomology of a modified $\lambda$-differential  3-Lie algebra with
coefficients in its representation.

Recall from \cite{Takhtajan} that
let $(\mathfrak{M};  \rho)$  be a representation of a   3-Lie algebra $(\mathfrak{A}, [-, -, -])$.
Denote the $n-$cochains   of $\mathfrak{A}$ with coefficients in representation $(\mathfrak{M};  \rho)$   by
\begin{align*}
\mathcal{C}_{\mathrm{3Lie}}^{n}(\mathfrak{A},\mathfrak{M}):=\mathrm{Hom}((\wedge^2\mathfrak{A})^{\otimes n-1}\wedge\mathfrak{A},\mathfrak{M}),n\geq 1.
\end{align*}

The   coboundary operator $\delta: \mathcal{C}_{\mathrm{3Lie}}^{n}(\mathfrak{A},\mathfrak{M})\rightarrow \mathcal{C}_{\mathrm{3Lie}}^{n+1}(\mathfrak{A},\mathfrak{M})$,  for $A_i=a_i\wedge b_{i}\in \wedge^2\mathfrak{A},a_{n+1}\in \mathfrak{A}$ and $f\in \mathcal{C}_{\mathrm{3Lie}}^{n}(\mathfrak{A},\mathfrak{M})$, as
\begin{align*}
&\delta f(A_1,\cdots, A_{n},a_{n+1})\\
=&(-1)^{n+1}\big(\rho(b_n,a_{n+1})f(A_1,\cdots, A_{n-1},a_{n})+\rho(a_{n+1},a_n)f(A_1,\cdots, A_{n-1},b_{n})\big)\\
&+\sum_{i=1}^n(-1)^{i+1}\rho(a_i,b_i)f(A_1,\cdots, A_{i-1},A_{i+1},\cdots, A_{n},a_{n+1})\\
&+\sum_{i=1}^n(-1)^{i}f(A_1,\cdots, A_{i-1},A_{i+1},\cdots, A_{n},[a_i,b_i,a_{n+1}])\\
&+\sum_{1\leq i<k\leq n }(-1)^{i}f(A_1,\cdots, A_{i-1},A_{i+1},\cdots, A_{k-1},[a_i,b_i,a_{k}]\wedge b_k+a_k\wedge [a_i,b_i,b_{k}],\cdots, A_{n},a_{n+1}),
\end{align*}
it was proved that $\delta\circ\delta=0.$

\begin{lemma}\label{lemma:cochain map}
Let $(\mathfrak{M};  \rho, \mathrm{d}_\mathfrak{M})$  be a representation of a modified $\lambda$-differential  3-Lie algebra $(\mathfrak{A}, [-, -, -],\mathrm{d})$.
For any $n\geq 1$, we define a linear map $\Phi: \mathcal{C}_{\mathrm{3Lie}}^{n}(\mathfrak{A},\mathfrak{M})\rightarrow \mathcal{C}_{\mathrm{3Lie}}^{n}(\mathfrak{A},\mathfrak{M})$ by
\begin{align*}
\Phi f(A_1,\cdots, A_{n-1},a_{n})=&\sum_{i=1}^{n-1}f(A_1,\cdots, A_{i-1},  \mathrm{d}(a_i)\wedge b_i+a_i\wedge \mathrm{d}(b_i),A_{i+1},\cdots, A_{n-1},a_{n})\\
&+f(A_1,\cdots, A_{n-1},\mathrm{d}(a_{n}))+(n-1)\lambda f(A_1,\cdots, A_{n-1},a_{n})\\
&-\mathrm{d}_\mathfrak{M}(f(A_1,\cdots, A_{n-1},a_{n})),
\end{align*}
for  any  $f\in \mathcal{C}_{\mathrm{3Lie}}^{n}(\mathfrak{A},\mathfrak{M})$ and $A_i=a_i\wedge b_{i}\in \wedge^2\mathfrak{A},i=1,\cdots,n-1, a_{n}\in \mathfrak{A}$.
Then $\Phi$ is a cochain map, i.e., $\Phi\circ\delta=\delta\circ\Phi.$
\end{lemma}

\begin{proof}
It follows by a straightforward tedious calculations.
\end{proof}

Let $(\mathfrak{M};  \rho, \mathrm{d}_\mathfrak{M})$ be a representation of the modified $\lambda$-differential 3-Lie algebra $(\mathfrak{A}, [-, -, -],\mathrm{d})$,
we define $n$-cochains for modified $\lambda$-differential  3-Lie algebra as follows:
\begin{equation*}\label{eq:dac}
\mathcal{C}_{\mathrm{md3Lie_\lambda}}^{n}(\mathfrak{A},\mathfrak{M}):=
\begin{cases}
\mathcal{C}^{n}_{\mathrm{3Lie}}(\mathfrak{A},\mathfrak{M})\oplus \mathcal{C}^{n-1}_{\mathrm{3Lie}}(\mathfrak{A},\mathfrak{M}),&n\geq 2,\\
\mathrm{Hom}(\mathfrak{A},\mathfrak{M}),&n=1.
\end{cases}
\end{equation*}

We define a  linear map  $\partial:\mathcal{C}_{\mathrm{md3Lie_\lambda}}^{n}(\mathfrak{A},\mathfrak{M})\rightarrow \mathcal{C}_{\mathrm{md3Lie_\lambda}}^{n+1}(\mathfrak{A},\mathfrak{M})$  by
\begin{align*}
\partial(f)&=(\delta f, -\Phi f), ~~~~~~~~~ ~~\mathrm{if}~~ f\in \mathcal{C}_{\mathrm{md3Lie_\lambda}}^{1}(\mathfrak{A},\mathfrak{M});\\
\partial(f,g)&=(\delta f, \delta g+(-1)^n \Phi f), ~ \mathrm{if} ~~(f,g)\in \mathcal{C}_{\mathrm{md3Lie_\lambda}}^{n}(\mathfrak{A},\mathfrak{M}).
\end{align*}

In view of Lemma \ref{lemma:cochain map}, we have the following theorem.

\begin{theorem}\label{thm: cochain complex for differential algebras}
The linear map $\partial $ is a coboundary operator, that is, $\partial\circ\partial=0.$
\end{theorem}

Therefore, we obtain a   cochain complex $(\mathcal{C}_{\mathrm{md3Lie_\lambda}}^*(\mathfrak{A}, \mathfrak{M}),\partial)$,
for $n\geq 2$,  we denote the set of $n$-cocycles by
$\mathcal{Z}^{n}_{\mathrm{md3Lie_\lambda}}(\mathfrak{A},\mathfrak{M})=\big\{(f,g)\in\mathcal{C}^{n}_{\mathrm{md3Lie_\lambda}}(\mathfrak{A},\mathfrak{ M})~|~$
$\partial (f,g)=0 \big\}$, the set of $n$-coboundaries by $\mathcal{B}^{n}_{\mathrm{md3Lie_\lambda}}(\mathfrak{A},\mathfrak{M})=\{\partial (f,g)~|~(f,g)\in\mathcal{C}^{n-1}_{\mathrm{md3Lie_\lambda}}(\mathfrak{A},\mathfrak{ M})\}$ and the $n$-th cohomology group of the  modified $\lambda$-differential  3-Lie algebra  $(\mathfrak{A}, [-, -, -],\mathrm{d})$  with coefficients in the representation    $(\mathfrak{M};  \rho, \mathrm{d}_\mathfrak{M})$ by $\mathcal{H}^{n}_{\mathrm{md3Lie_\lambda}}(\mathfrak{A},\mathfrak{M})=\mathcal{Z}^{n}_{\mathrm{md3Lie_\lambda}}(\mathfrak{A},\mathfrak{ M})/\mathcal{B}^{n}_{\mathrm{md3Lie_\lambda}}(\mathfrak{A},\mathfrak{M})$.

Lastly, we calculate the 1-cocycle and 2-cocycle.

For $f\in \mathcal{C}_{\mathrm{md3Lie_\lambda}}^{1}(\mathfrak{A},\mathfrak{M})$, $f$ is a 1-cocycle   if $\partial(f)=(\delta f, -\Phi f)=0,$ i.e.,
$$\rho(b_1,a_2)f(a_1)+\rho(a_2,a_1)f(b_1)+\rho(a_1,b_1)f(a_2)-f([a_1,b_1,a_2])=0$$
and
$$\mathrm{d}_\mathfrak{M}(f(a_1))-f(\mathrm{d}(a_1))=0.$$

For $(f,g)\in \mathcal{C}_{\mathrm{md3Lie_\lambda}}^{2}(\mathfrak{A},\mathfrak{M})$, $(f,g)$ is a 2-cocycle  if $\partial(f,g)=(\delta f, \delta g+\Phi f)=0,$ i.e.,
\begin{align*}
&-\rho(b_2,a_3)f(a_1,b_1,a_2)-\rho(a_3,a_2)f(a_1,b_1,b_2)+\rho(a_1,b_1)f(a_2,b_2,a_3)-\rho(a_2,b_2)f(a_1,b_1,a_3)\\
&-f(a_2,b_2,[a_1,b_1,a_3])+f(a_1,b_1,[a_2,b_2,a_3])-f([a_1,b_1,a_2],b_2,a_3)-f(a_2,[a_1,b_1,b_2],a_3)=0
\end{align*}
and
\begin{align*}
&\rho(b_1,a_2)f(a_1)+\rho(a_2,a_1)f(b_1)+\rho(a_1,b_1)f(a_2)-f([a_1,b_1,a_2])\\
&+f(\mathrm{d}(a_1),b_1,a_2)+f(a_1,\mathrm{d}(b_1),a_2)+f(a_1,b_1,d(a_2))+\lambda f(a_1,b_1,a_2)-\mathrm{d}_\mathfrak{M}(f(a_1,b_1,a_2))=0.
\end{align*}

\section{Deformations of modified $\lambda$-differential 3-Lie algebras}\label{sec:def}
\def\theequation{\arabic{section}.\arabic{equation}}
\setcounter{equation} {0}

 In this section, we  consider linear deformations  of  the modified $\lambda$-differential 3-Lie algebra.

 Let $(\mathfrak{A}, [-, -, -],\mathrm{d})$ be a   modified $\lambda$-differential 3-Lie algebra. Denote  $\nu_0=[-, -, -]$ and $\mathrm{d}_0=\mathrm{d}$.
Consider a family of linear maps:
$$\nu_t=\nu_0+t\nu_1+t^2\nu_2, \, \, \nu_1,\nu_2\in \mathcal{C}^2_{\mathrm{3Lie}}(\mathfrak{A}, \mathfrak{A}),\quad
 \mathrm{d}_t=\mathrm{d}_0+t   \mathrm{d}_1, \, \, \mathrm{d}_1\in \mathcal{C}^1_{\mathrm{3Lie}}(\mathfrak{A}, \mathfrak{A}).$$

\begin{defn}
A  linear deformation of the modified $\lambda$-differential 3-Lie algebra $(\mathfrak{A}, [-, -, -],\mathrm{d})$  is a pair $(\nu_t, \mathrm{d}_t)$ which endows $(\mathfrak{A}[[t]], \nu_t, \mathrm{d}_t)$ with the modified $\lambda$-differential  3-Lie algebra.
\end{defn}


 \begin{prop}\label{prop:fddco}

 The pair $(\nu_t, \mathrm{d}_t)$ generates a linear deformation    of the modified $\lambda$-differential  3-Lie algebra $(\mathfrak{A}, [-, -, -],\mathrm{d})$   if and only if the following equations hold:
\begin{align}
\sum_{i+j=n}\nu_i(a_1,a_2,\nu_{j}(a_3,a_4,a_5))=&\sum_{i+j=n}\nu_i(\nu_{j}(a_1,a_2,a_3),a_4,a_5)+\sum_{i+j=n}\nu_i(a_3,\nu_{j}(a_1,a_2,a_4),a_5)\nonumber\\
&+\sum_{i+j=n}\nu_i(a_3,a_4,\nu_{j}(a_1,a_2,a_5)),\label{3.1}\\
\sum_{i+l=n}\mathrm{d}_l(\nu_{i}(a_1,a_2,a_3))=&\sum_{i+l=n}\nu_{i}(\mathrm{d}_l(a_1),a_2,a_3)+\sum_{i+l=n}\nu_{i}(a_1,\mathrm{d}_l(a_2),a_3)\nonumber\\
 &+\sum_{i+l=n}\nu_{i}(a_1,a_2,\mathrm{d}_l(a_3))+\lambda \nu_n(a_1,a_2,a_3),\label{3.2}
\end{align}
for any  $a_1,a_2,a_3,a_4,a_5\in \mathfrak{A}$ and $i,j=0,1,2,l=0,1$.
\end{prop}

\begin{proof}
$(\mathfrak{A}[[t]], \nu_t, \mathrm{d}_t)$ is a modified $\lambda$-differential  3-Lie algebra if and only if
\begin{align}
&\nu_t(a_1,a_2,\nu_t(a_3,a_4,a_5))\nonumber\\
&=\nu_t(\nu_t(a_1,a_2,a_3),a_4,a_5)+\nu_t(a_3,\nu_t(a_1,a_2,a_4),a_5)+\nu_t(a_3,a_4,\nu_t(a_1,a_2,a_5)),\label{3.3}\\
 &\mathrm{d}_t(\nu_t(a_1,a_2,a_3))\nonumber\\
 &=\nu_t(\mathrm{d}_t(a_1),a_2,a_3)+\nu_t(a_1,\mathrm{d}_t(a_2),a_3)+\nu_t(a_1,a_2,\mathrm{d}_t(a_3))+\lambda \nu_t(a_1,a_2,a_3).\label{3.4}
\end{align}
Comparing the coefficients of $t^n$ on both sides of the above equations, Eqs. \eqref{3.3} and \eqref{3.4} are equivalent to  Eqs. \eqref{3.1} and \eqref{3.2} respectively.
\end{proof}

 \begin{coro}\label{prop:fddco}
Let $(\mathfrak{A}[[t]], \nu_t, \mathrm{d}_t)$ be a linear deformation   of the modified $\lambda$-differential 3-Lie algebra $(\mathfrak{A}, \nu_0,\mathrm{d})$.  Then $(\nu_1, \mathrm{d}_1)$ is a 2-cocycle of $(\mathfrak{A}, [-, -, -],\mathrm{d})$ with the coefficient  in the adjoint representation $(\mathfrak{A}; ad,\mathrm{ d})$.	
\end{coro}
\begin{proof}
  For $n =1$, Eqs.~\eqref{3.1} and \eqref{3.2} are equivalent to
\begin{align*}
&\nu_1(a_1,a_2,[a_3,a_4,a_5])+[a_1,a_2,\nu_1(a_3,a_4,a_5)]\\
 =& \nu_1([a_1,a_2,a_3],a_4,a_5)+[\nu_1(a_1,a_2,a_3),a_4,a_5]+\nu_1(a_3,[a_1,a_2,a_4],a_5)+[a_3,\nu_1(a_1,a_2,a_4),a_5]\\
 &+\nu_1(a_3,a_4,[a_1,a_2,a_5])+[a_3,a_4,\nu_1(a_1,a_2,a_5)],\\
&\mathrm{d}_1([a_1,a_2,a_3])+ \mathrm{d}(\nu_1(a_1,a_2,a_3))\\
 =& [\mathrm{d}_1(a_1),a_2,a_3]+\nu_1(\mathrm{d}(a_1),a_2,a_3)+[a_1,\mathrm{d}_1(a_2),a_3]+\nu_1(a_1,\mathrm{d}(a_2),a_3)+[a_1,a_2,\mathrm{d}_1(a_3)]\\
 &+\nu_1(a_1,a_2,\mathrm{d}(a_3))+\lambda\nu_1(a_1,a_2,a_3),
 \end{align*}
which imply that  $\delta \nu_1=0,\delta \mathrm{d}_1+ \Phi\nu_1=0$ respectively.  Hence, 	$(\nu_1, \mathrm{d}_1)$ is a 2-cocycle of $(\mathfrak{A}, [-, -, -],\mathrm{d})$ with the coefficient  in the adjoint representation $(\mathfrak{A};ad,\mathrm{d})$.	
\end{proof}

\begin{defn}
The $2$-cocycle $(\nu_1,\mathrm{d}_1)$ is called the   infinitesimal  of the linear deformation $(\mathfrak{A}[[t]],\nu_t,\mathrm{d}_t)$ of $(\mathfrak{A}, [-, -, -],\mathrm{d})$.
\end{defn}

\begin{defn}
(i)
Two linear deformations $(\mathfrak{A}[[t]],\nu_t,\mathrm{d}_t)$ and $(\mathfrak{A}[[t]],\nu'_t,\mathrm{d}'_t)$ of the modified $\lambda$-differential 3-Lie algebra $(\mathfrak{A}, [-, -, -],\mathrm{d})$ are said to be   equivalent  if there exists a linear map   $N:\mathfrak{A}\rightarrow\mathfrak{A}$, such that $N_t=\mathrm{id}_\mathfrak{A}+tN$ satisfying
\begin{align}
N_t(\mathrm{d}_t(a_1))=&\mathrm{d}'_t(N_t(a_1)) ,\label{3.5}\\
N_t\nu_t(a_1,a_2,a_3)=&\nu'_t(N_t(a_1),N_t(a_2),N_t(a_3)),\label{3.6}
\end{align}
for any $a_1,a_2,a_3\in \mathfrak{A}$.\\
(ii) A linear deformation  $(\mathfrak{A}[[t]],\nu_t,\mathrm{d}_t)$ of the modified $\lambda$-differential 3-Lie algebra $(\mathfrak{A}, [-, -, -],\mathrm{d})$ is  said to be   trivial  if  $(\mathfrak{A}[[t]],\nu_t,\mathrm{d}_t)$ is equivalent to $(\mathfrak{A}, [-, -, -],\mathrm{d})$.
\end{defn}

Comparing the coefficients of $t$ on both sides of the above Eqs.  \eqref{3.5} and \eqref{3.6}, we have
\begin{align*}
\nu_1(a_1,a_2,a_3)-\nu'_1(a_1,a_2,a_3)&=[Na_1,a_2,a_3]+[a_1,Na_2,a_3]+[a_1,a_2,Na_3]-N[a_1,a_2,a_3],\\
\mathrm{d}_1(a)-\mathrm{d}'_1(a)&=\mathrm{d}(N_1a)-N_1\mathrm{d}(a).
\end{align*}
Thus, we have the following theorem.
\begin{theorem}
The infinitesimals of two equivalent linear  deformations of $(\mathfrak{A}, [-, -, -],\mathrm{d})$ are in the same  cohomological class in $\mathcal{H}_{\mathrm{md3Lie_{\lambda}}}^2(\mathfrak{A},\mathfrak{A})$.
\end{theorem}

Let  $(\mathfrak{A}[[t]],\nu_t,\mathrm{d}_t)$  be a trivial deformation of $(\mathfrak{A}, [-, -, -],\mathrm{d})$. Then there exists a linear map   $N:\mathfrak{A}\rightarrow\mathfrak{A}$, such that $N_t=\mathrm{id}_\mathfrak{A}+tN$ satisfying
\begin{align}
N_t(\mathrm{d}(a_1))=&\mathrm{d}(N_t(a_1)) ,\label{3.7}\\
N_t\nu_t(a_1,a_2,a_3)=&[N_t(a_1),N_t(a_2),N_t(a_3)].\label{3.8}
\end{align}
Compare the coefficients of $t^i (1\leq i\leq 3)$ on both sides of Eqs. \eqref{3.7} and \eqref{3.8}, and we can get
\begin{align}
N \mathrm{d}(a_1) =&\mathrm{d} (Na_1),\label{3.9}\\
\nu_1(a_1,a_2,a_3)+N[a_1,a_2,a_3]=&[Na_1,  a_2, a_3]+[ a_1, Na_2, a_3]+[ a_1,  a_2,Na_3],\label{3.10}\\
\nu_2(a_1,a_2,a_3)+N\nu_1(a_1,a_2,a_3)=&[Na_1,  Na_2, a_3]+[ a_1, Na_2, Na_3]+[ Na_1,  a_2,Na_3],\label{3.11}\\
N\nu_2(a_1,a_2,a_3)=&[Na_1,  Na_2, Na_3]\label{3.12}.
\end{align}

Thus, from a trivial deformation, we can get the following definition of Nijenhuis operator.

\begin{defn}
Let $(\mathfrak{A}, [-, -, -],\mathrm{d})$ be a   modified $\lambda$-differential 3-Lie algebra. A linear map
$N: \mathfrak{A}\rightarrow\mathfrak{A}$ is called a Nijenhuis operator if the following equations hold:
\begin{align}
N\circ\mathrm{d}=&\mathrm{d}\circ N, \label{3.13}\\
[Na_1,Na_2,Na_3]=&N([a_1,Na_2,Na_3]+[Na_1,a_2,Na_3]+[Na_1,Na_2,a_3])\nonumber\\
&-N^2([Na_1,a_2,a_3]+[a_1,Na_2,a_3]+[a_1,a_2,Na_3])+N^3[a_1,a_2,a_3],\label{3.14}
\end{align}
for any $a_1,a_2,a_3\in\mathfrak{A}.$
\end{defn}

 \begin{prop}
Let $(\mathfrak{A}, [-, -, -],\mathrm{d})$ be a   modified $\lambda$-differential 3-Lie algebra, and $N: \mathfrak{A}\rightarrow\mathfrak{A}$   a Nijenhuis operator. Then $(\mathfrak{A}, [-, -, -]_N,\mathrm{d})$  is  a   modified $\lambda$-differential 3-Lie algebra,
where
\begin{align*}
[a_1,a_2,a_3]_N=&[a_1,Na_2,Na_3]+[Na_1,a_2,Na_3]+[Na_1,Na_2,a_3]\\
&-N([Na_1,a_2,a_3]+[a_1,Na_2,a_3]+[a_1,a_2,Na_3])+N^2[a_1,a_2,a_3].
\end{align*}
\end{prop}

\begin{proof}
In the light of \cite{Liu}, $(\mathfrak{A}, [-, -, -]_N)$  is  a    3-Lie algebra.  Next we prove that  $\mathrm{d}$ is a modified $\lambda$-differential  operator of $(\mathfrak{A}, [-, -, -]_N)$,
for any $a_1,a_2,a_3\in\mathfrak{A},$ by Eqs. \eqref{2.2} and \eqref{3.13}, we have
\begin{align*}
&\mathrm{d}[a_1,a_2,a_3]_N\\
=&\mathrm{d}[a_1,Na_2,Na_3]+\mathrm{d}[Na_1,a_2,Na_3]+\mathrm{d}[Na_1,Na_2,a_3]\\
&-N(\mathrm{d}[Na_1,a_2,a_3]+\mathrm{d}[a_1,Na_2,a_3]+\mathrm{d}[a_1,a_2,Na_3])+N^2\mathrm{d}[a_1,a_2,a_3]\\
=&[\mathrm{d}(a_1), Na_2,  Na_3]+[a_1,  N\mathrm{d}(a_2),  Na_3]+[a_1, Na_2,  N\mathrm{d}(a_3)]+\lambda[a_1, Na_2,  Na_3]\\
&+[N\mathrm{d}(a_1), a_2, N a_3]+[Na_1,  \mathrm{d}(a_2),  Na_3]+[Na_1, a_2,  N\mathrm{d}(a_3)]+\lambda[Na_1, a_2,  Na_3]\\
&+[N\mathrm{d}(a_1), Na_2,  a_3]+[Na_1,  N\mathrm{d}(a_2),  a_3]+[Na_1, Na_2,  \mathrm{d}(a_3)]+\lambda[Na_1, Na_2,  a_3]\\
&-N([N\mathrm{d}(a_1), a_2,  a_3]+[Na_1,  \mathrm{d}(a_2),  a_3]+[Na_1, a_2,  \mathrm{d}(a_3)]+\lambda[Na_1, a_2,  a_3])\\
&-N([\mathrm{d}(a_1), Na_2,  a_3]+[a_1,  N\mathrm{d}(a_2),  a_3]+[a_1, Na_2,  \mathrm{d}(a_3)]+\lambda[a_1, Na_2,  a_3])\\
&-N([\mathrm{d}(a_1), a_2,  Na_3]+[a_1,  \mathrm{d}(a_2),  Na_3]+[a_1, a_2,  N\mathrm{d}(a_3)]+\lambda[a_1, a_2,  Na_3])\\
&+N^2([\mathrm{d}(a_1), a_2,  a_3]+[a_1,  \mathrm{d}(a_2),  a_3]+[a_1, a_2,  \mathrm{d}(a_3)]+\lambda[a_1, a_2,  a_3])\\
=&[\mathrm{d}(a_1), a_2,  a_3]_N+[a_1,  \mathrm{d}(a_2),  a_3]_N+[a_1, a_2,  \mathrm{d}(a_3)]_N+\lambda[a_1, a_2,  a_3]_N.
\end{align*}
So we get the conclusion.
\end{proof}

\begin{defn}
A linear map $R: \mathfrak{M}\rightarrow\mathfrak{A}$ is called an $\mathcal{O}$-operator on the modified $\lambda$-differential 3-Lie algebra $(\mathfrak{A}, [-, -, -],\mathrm{d})$
with respect to the representation $(\mathfrak{M};  \rho, \mathrm{d}_\mathfrak{M})$    if the following equations hold:
\begin{align*}
R\circ\mathrm{d}_\mathfrak{M}=&\mathrm{d}\circ R, \\
[Rv_1,Rv_2,Rv_3]=&R(\rho(Rv_1,Rv_2)v_3+\rho(Rv_2,Rv_3)v_1+\rho(Rv_3, Rv_1)v_2),
\end{align*}
for any $v_1,v_2,v_3\in\mathfrak{M}.$
\end{defn}

 \begin{remark}
Obviously, An invertible linear map $R: \mathfrak{M}\rightarrow\mathfrak{A}$ is  an $\mathcal{O}$-operator if and only if $R^{-1}$ is a 1-cocycle of the  modified $\lambda$-differential  3-Lie algebra  $(\mathfrak{A}, [-, -, -],\mathrm{d})$  with coefficients in the representation    $(\mathfrak{M};  \rho, \mathrm{d}_\mathfrak{M})$.
\end{remark}

 \begin{prop}
Let $(\mathfrak{M};  \rho, \mathrm{d}_\mathfrak{M})$  be a  representation  of a   modified $\lambda$-differential 3-Lie algebra $(\mathfrak{A}, [-, -, -],\mathrm{d})$.   Then $R: \mathfrak{M}\rightarrow\mathfrak{A}$ is  an  $\mathcal{O}$-operator if and only if $\overline{R}=\left(
        \begin{array}{cc}
          0 & R \\
          0 & 0 \\
        \end{array}
      \right):\mathfrak{A}\oplus\mathfrak{M}\rightarrow \mathfrak{A}\oplus\mathfrak{M}$ is a Nijenhuis operator on semidirect product  modified $\lambda$-differential 3-Lie algebra  $(\mathfrak{A} \oplus \mathfrak{M},[-, -, -]_{\rho},\mathrm{d} \oplus \mathrm{d}_\mathfrak{M})$.
\end{prop}

\begin{proof}
For any $a_1, a_2, a_3\in \mathfrak{A}$ and $u_1, u_2, u_3\in \mathfrak{M}$, by $\overline{R}^2=0,$ we have
\begin{align*}
&(\mathrm{d} \oplus \mathrm{d}_\mathfrak{M})\overline{R}(a_1+u_1)=(\mathrm{d} \oplus \mathrm{d}_\mathfrak{M})(Ru_1+0)=\mathrm{d}(Ru_1)+0,\\
&\overline{R}(\mathrm{d} \oplus \mathrm{d}_\mathfrak{M})(a_1+u_1)=\overline{R}(\mathrm{d}(a_1)+\mathrm{d}_\mathfrak{M}(u_1))=R\mathrm{d}_\mathfrak{M}(u_1)+0,\\
&\overline{R}([a_1+u_1, \overline{R}(a_2+u_2), \overline{R}(a_3+u_3)]_{\rho}+[\overline{R}(a_1+u_1), a_2+u_2, \overline{R}(a_3+u_3)]_{\rho}\\
&+[\overline{R}(a_1+u_1), \overline{R}(a_2+u_2), a_3+u_3]_{\rho})-[\overline{R}(a_1+u_1), \overline{R}(a_2+u_2), \overline{R}(a_3+u_3)]_{\rho}\\
=&\overline{R}([a_1+u_1, Ru_2+0, Ru_3+0]_{\rho}+[Ru_1+0, a_2+u_2, Ru_3+0]_{\rho}+[Ru_1+0, Ru_2+0, a_3+u_3]_{\rho})\\
&-[Ru_1+0, Ru_2+0, Ru_3+0]_{\rho}\\
=&\overline{R}([a_1, Ru_2, Ru_3]+\rho(Ru_2, Ru_3)u_1+
[Ru_1, a_2, Ru_3]+\rho(Ru_3, Ru_1)u_2+[Ru_1, Ru_2 a_3]+\rho(Ru_1, Ru_2)u_3)\\
&-[Ru_1, Ru_2, Ru_3]+0\\
=&R(\rho(Ru_2, Ru_3)u_1+\rho(Ru_3, Ru_1)u_2+\rho(Ru_1, Ru_2)u_3)-[Ru_1, Ru_2, Ru_3]+0,
\end{align*}
which implies that $R$ is  an  $\mathcal{O}$-operator if and only if $\overline{R}$ is a Nijenhuis operator.
\end{proof}

\section{Abelian extensions of modified $\lambda$-differential 3-Lie algebras} \label{sec:aext}
\def\theequation{\arabic{section}.\arabic{equation}}
\setcounter{equation} {0}
In this section, we study abelian extensions of modified $\lambda$-differential 3-Lie algebras and show that they are classified   by the second cohomology groups.

Notice that a vector space $\mathfrak{M}$ together with a linear map $\mathrm{d}_\mathfrak{M}$ is naturally a
modified $\lambda$-differential  3-Lie algebra where  the bracket on $\mathfrak{M}$ is defined to be  $[-,-,-]_\mathfrak{M}=0.$

\begin{defn}
An abelian extension of $(\mathfrak{A}, [-,-,-], \mathrm{d})$ by  $(\mathfrak{M},[-,-,-]_\mathfrak{M}, \mathrm{d}_\mathfrak{M})$
 is a short exact sequence of homomorphisms of modified $\lambda$-differential  3-Lie algebras
  $$0\rightarrow  (\mathfrak{M},[-,-,-]_\mathfrak{M}, \mathrm{d}_\mathfrak{M})\stackrel{i}{\longrightarrow} (\widehat{\mathfrak{A}}, [-,-,-]_{\widehat{\mathfrak{A}}}, \widehat{\mathrm{d}})\stackrel{p}{\longrightarrow} (\mathfrak{A}, [-,-,-], \mathrm{d})\rightarrow 0,$$
i.e., there exists a commutative diagram:
$$\begin{CD}
0@>>> {\mathfrak{M}} @>i >> \widehat{\mathfrak{A}} @>p >> \mathfrak{A} @>>>0\\
@. @V {\mathrm{d}_\mathfrak{M}} VV @V \widehat{\mathrm{d}} VV @V \mathrm{d} VV @.\\
0@>>> {\mathfrak{M}} @>i >> \widehat{\mathfrak{A}} @>p >> \mathfrak{A} @>>>0,
\end{CD}$$
where the modified $\lambda$-differential  3-Lie algebras $(\widehat{\mathfrak{A}}, [-,-,-]_{\widehat{\mathfrak{A}}}, \widehat{\mathrm{d}})$  satisfies   $[-,u, v]_{\widehat{\mathfrak{A}}}=0$, for all $u,v\in \mathfrak{M}$.

We will call $(\widehat{\mathfrak{A}}, [-,-,-]_{\widehat{\mathfrak{A}}}, \widehat{\mathrm{d}})$  an abelian extension of $(\mathfrak{A}, [-,-,-], \mathrm{d})$ by $(\mathfrak{M},[-,-,-]_\mathfrak{M}, \mathrm{d}_\mathfrak{M})$.
\end{defn}

A   section  of an abelian extension $(\widehat{\mathfrak{A}}, [-,-,-]_{\widehat {\mathfrak{A}}}, \widehat {\mathrm{d}})$ of $(\mathfrak{A}, [-,-,-], \mathrm{d})$ by  $(\mathfrak{M},[-,-,-]_\mathfrak{M},\mathrm{d}_\mathfrak{M})$ is a linear map $s:\mathfrak{A}\rightarrow \widehat{\mathfrak{A}}$ such that $p\circ s=\mathrm{id}_\mathfrak{A}$.

Let $(\mathfrak{M};  \rho, \mathrm{d}_\mathfrak{M})$  be a  representation  of a   modified $\lambda$-differential 3-Lie algebra $(\mathfrak{A}, [-, -, -],\mathrm{d})$.
Assume  that $(f,g)\in\mathcal{C}_{\mathrm{md3Lie_\lambda}}^2(\mathfrak{A}, \mathfrak{M})$. Define $[-,-,-]_{\rho f}: \wedge^3 (\mathfrak{A}\oplus\mathfrak{M})\rightarrow \mathfrak{A}\oplus\mathfrak{M}$
and   $ \mathrm{d}_{g}: \mathfrak{A}\oplus\mathfrak{M}\rightarrow \mathfrak{A}\oplus\mathfrak{M}$ respectively by
\begin{align}
&[a_1+u_1, a_2+u_2,a_3+u_3]_{\rho f} \nonumber\\
&=[a_1, a_2, a_3]+\rho(a_2,a_3)u_1+\rho(a_3,a_1)u_2+\rho(a_1,a_2)u_3+f(a_1, a_2, a_3), \label{4.1} \\
&\mathrm{d}_g(a_1+u_1)=\mathrm{d}(a_1)+\mathrm{d}_\mathfrak{M}(u_1)+g(a_1),\,\forall a_1, a_2, a_3\in \mathfrak{A},\ u_1, u_2, u_3\in \mathfrak{M}.\label{4.2}
\end{align}

\begin{prop} \label{prop:2-cocycle}
The triple $(\mathfrak{A}\oplus \mathfrak{M},[-,-,-]_{\rho f},\mathrm{d}_{g})$ is a modified $\lambda$-differential 3-Lie algebra  if and only if
$(f,g)$ is a 2-cocycle   in the cohomology of the modified $\lambda$-differential 3-Lie algebra $(\mathfrak{A}, [-, -, -],\mathrm{d})$ with the coefficient  in $(\mathfrak{M};  \rho, \mathrm{d}_\mathfrak{M})$.
 In this case, $$0\rightarrow  (\mathfrak{M},[-,-,-]_\mathfrak{M}, \mathrm{d}_\mathfrak{M})\hookrightarrow (\mathfrak{A}\oplus \mathfrak{M},[-,-,-]_{\rho f},\mathrm{d}_{g})\stackrel{p}{\longrightarrow} (\mathfrak{A}, [-,-,-], \mathrm{d})\rightarrow 0$$ is an abelian extension.
\end{prop}

\begin{proof}
 $(\mathfrak{A}\oplus \mathfrak{M},[-,-,-]_{\rho f},\mathrm{d}_{g})$ is a modified $\lambda$-differential 3-Lie algebra if and only if
\begin{align}
&[[a_1+u_1, a_2+u_2,a_3+u_3]_{\rho f}, a_4+u_4,a_5+u_5]_{\rho f}+[a_3+u_3,[a_1+u_1, a_2+u_2, a_4+u_4]_{\rho f},a_5+u_5]_{\rho f}\nonumber\\
&+[ a_3+u_3, a_4+u_4, [a_1+u_1, a_2+u_2,a_5+u_5]_{\rho f}]_{\rho f}- [a_1+u_1, a_2+u_2, [a_3+u_3, a_4+u_4,a_5+u_5]_{\rho f}]_{\rho f} \nonumber\\
&=0, \label{4.3}\\
&[\mathrm{d}_g(a_1+u_1), a_2+u_2,a_3+u_3]_{\rho f}+[a_1+u_1, \mathrm{d}_g (a_2+u_2),a_3+u_3]_{\rho f}+[a_1+u_1, a_2+u_2,\mathrm{d}_g (a_3+u_3)]_{\rho f}\nonumber\\
&+\lambda[a_1+u_1, a_2+u_2,a_3+u_3]_{\rho f}-\mathrm{d}_g[a_1+u_1, a_2+u_2,a_3+u_3]_{\rho f}\nonumber\\
&=0, \label{4.4}
\end{align}
for any $a_1, a_2, a_3,a_4, a_5\in \mathfrak{A}, \ u_1, u_2, u_3,\ u_4, u_5\in \mathfrak{M}$.
Furthermore, Eqs. \eqref{4.3} and \eqref{4.4} are equivalent to
\begin{align}
&\rho(a_4,a_5)f(a_1,a_2,a_3)+f( [a_1,a_2,a_3],a_4,a_5)+\rho(a_5,a_3)f(a_1,a_2,a_4)+f(a_3, [a_1,a_2,a_4],a_5)\nonumber\\
&+\rho(a_3,a_4)f(a_1,a_2,a_5)+f(a_3,a_4, [a_1,a_2,a_5])-\rho(a_1,a_2)f(a_3,a_4,a_5)-f(a_1,a_2, [a_3,a_4,a_5])\nonumber\\
&=0 \label{4.5}\\
&\rho(a_2,a_3)g(a_1)+f( \mathrm{d}(a_1),a_2,a_3)+ \rho(a_3,a_1)g(a_2)+f(a_1,\mathrm{d}(a_2),a_3)+\rho(a_1,a_2)g(a_3)\nonumber\\
&+f(a_1,a_2, \mathrm{d}(a_3))+\lambda f(a_1,a_2,a_3)-\mathrm{d}_\mathfrak{M}(f(a_1,a_2,a_3))-g([a_1,a_2,a_3])=0, \label{4.6}
\end{align}
using Eqs. \eqref{4.5} and  \eqref{4.6}, we get $\delta f=0$ and $\delta g+\Phi f=0$, respectively.
Therefore, $\partial(f,g)=(\delta f,\delta g+\Phi f)=0,$  which implies that $(f,g)$ is a  2-cocycle.

Conversely, if $(f,g)$ satisfying Eqs.~\eqref{4.5} and \eqref{4.6}, then  $(\mathfrak{A}\oplus \mathfrak{M},[-,-,-]_{\rho f},\mathrm{d}_{g})$ is a modified $\lambda$-differential 3-Lie algebra.
\end{proof}

  Let $(\widehat{\mathfrak{A}}, [-,-,-]_{\widehat {\mathfrak{A}}}, \widehat {\mathrm{d}})$ be an abelian extension of $(\mathfrak{A}, [-,-,-], \mathrm{d})$ by  $(\mathfrak{M},[-,-,-]_\mathfrak{M}, \mathrm{d}_\mathfrak{M})$ and $s:\mathfrak{A}\rightarrow\widehat{\mathfrak{A}}$  a section. Define  $\varrho: \wedge^2\mathfrak{A} \rightarrow \mathrm{End}(\mathfrak{M})$, $\upsilon:\wedge^3\mathfrak{A}\rightarrow \mathfrak{M}$ and $\mu:\mathfrak{A}\rightarrow \mathfrak{M}$ respectively by
\begin{align*}
\varrho(a_1,a_2)u:&=[s(a_1),s(a_2),u]_{\widehat{\mathfrak{A}}},  \\
\upsilon(a_1,a_2,a_3):&=[s(a_1),s(a_2),s(a_3)]_{\widehat{\mathfrak{A}}}-s([a_1,a_2,a_3]),\\
\mu(a_1):&=\widehat{\mathrm{d}}(s(a_1))-s(\mathrm{d}(a_1)),\quad \forall a_1,a_2,a_3\in \mathfrak{A}, u\in \mathfrak{M}.
\end{align*}
Note that $\varrho$ is independent on the choice of $s$.
\begin{prop} \label{prop:r2-cocycle}
  With the above notations, $(\mathfrak{M},\varrho, \mathrm{d}_\mathfrak{M})$ is a representation over the modified $\lambda$-differential 3-Lie algebras   $(\mathfrak{A}, [-,-,-], \mathrm{d})$ and
$(\upsilon,\mu)$ is a 2-cocycle in the cohomology  of the modified $\lambda$-differential 3-Lie algebras   $(\mathfrak{A}, [-,-,-], \mathrm{d})$ with the coefficient  in $(\mathfrak{M}; \varrho, \mathrm{d}_\mathfrak{M})$.
Furthermore, the cohomological class of the 2-cocycle $[(\upsilon,\mu)]\in \mathcal{H}_{\mathrm{md3Lie_{\lambda}}}^2(\mathfrak{A},\mathfrak{M})$  is independent of the choice of sections of  $p$.
\end{prop}
\begin{proof}
First, for any  $a_1,a_2,a_3,a_4\in \mathfrak{A}, u\in \mathfrak{M}$, by Eq. \eqref{2.1},   we get
\begin{align*}
&\varrho(a_2,a_3)\varrho(a_1,a_4)u+\varrho(a_3,a_1)\varrho(a_2,a_4)u+\varrho(a_1,a_2)\varrho(a_3,a_4)u-\varrho([a_1,a_2,a_3],a_4)u\\
=&[s(a_2),s(a_3),[s(a_1),s(a_4),u]_{\widehat{\mathfrak{A}}}]_{\widehat{\mathfrak{A}}}+[s(a_3),s(a_1),[s(a_2),s(a_4),u]_{\widehat{\mathfrak{A}}}]_{\widehat{\mathfrak{A}}}\\
&+[s(a_1),s(a_2),[s(a_3),s(a_4),u]_{\widehat{\mathfrak{A}}}]_{\widehat{\mathfrak{A}}}-[s([a_1,a_2,a_3]),s(a_4),u]_{\widehat{\mathfrak{A}}}\\
=&[s(a_2),s(a_3),[s(a_1),s(a_4),u]_{\widehat{\mathfrak{A}}}]_{\widehat{\mathfrak{A}}}+[s(a_3),s(a_1),[s(a_2),s(a_4),u]_{\widehat{\mathfrak{A}}}]_{\widehat{\mathfrak{A}}}\\
&+[s(a_1),s(a_2),[s(a_3),s(a_4),u]_{\widehat{\mathfrak{A}}}]_{\widehat{\mathfrak{A}}}-[[s(a_1),s(a_2),s(a_3)]_{\widehat{\mathfrak{A}}}-\upsilon(a_1,a_2,a_3),s(a_4),u]_{\widehat{\mathfrak{A}}}\\
=&0,\\
&\varrho(a_3,a_4)\varrho(a_1,a_2)u+\varrho([a_1,a_2,a_3],a_4)u+\varrho(a_3,[a_1,a_2,a_4])u-\varrho(a_1,a_2)\varrho(a_3,a_4)u\\
=&[s(a_3),s(a_4),[s(a_1),s(a_2),u]_{\widehat{\mathfrak{A}}}]_{\widehat{\mathfrak{A}}}+[s([a_1,a_2,a_3]),s(a_4),u]_{\widehat{\mathfrak{A}}}+[s(a_3),s([a_1,a_2,a_4]),u]_{\widehat{\mathfrak{A}}}\\
&-[s(a_1),s(a_2),[s(a_3),s(a_4),u]_{\widehat{\mathfrak{A}}}]_{\widehat{\mathfrak{A}}}\\
=&[s(a_3),s(a_4),[s(a_1),s(a_2),u]_{\widehat{\mathfrak{A}}}]_{\widehat{\mathfrak{A}}}+[[s(a_1),s(a_2),s(a_3)]-\upsilon(a_1,a_2,a_3),s(a_4),u]_{\widehat{\mathfrak{A}}}\\
&+[s(a_3),[s(a_1),s(a_2),s(a_4)]-\upsilon(a_1,a_2,a_4),u]_{\widehat{\mathfrak{A}}}-[s(a_1),s(a_2),[s(a_3),s(a_4),u]_{\widehat{\mathfrak{A}}}]_{\widehat{\mathfrak{A}}}\\
=&0.
\end{align*}
In addition,  by Eq.   \eqref{2.2},   we have
\begin{align*}
&\mathrm{d}_\mathfrak{M}(\varrho(a_1,a_2)u)=\mathrm{d}_\mathfrak{M}([s(a_1),s(a_2),u]_{\widehat{\mathfrak{A}}})\\
=&[\widehat{\mathrm{d}}(s(a_1)),s(a_2),u]_{\widehat{\mathfrak{A}}}+[s(a_1),\widehat{\mathrm{d}}(s(a_2)),u]_{\widehat{\mathfrak{A}}}+[s(a_1),s(a_2),\mathrm{d}_\mathfrak{M}(u)]_{\widehat{\mathfrak{A}}}+\lambda [s(a_1),s(a_2),u]_{\widehat{\mathfrak{A}}} \\
=&[s(\mathrm{d}(a_1))+\mu(a_1),s(a_2),u]_{\widehat{\mathfrak{A}}}+[s(a_1),s(\mathrm{d}(a_2))+\mu(a_2),u]_{\widehat{\mathfrak{A}}}+[s(a_1),s(a_2),\mathrm{d}_\mathfrak{M}(u)]_{\widehat{\mathfrak{A}}}\\
&+\lambda [s(a_1),s(a_2),u]_{\widehat{\mathfrak{A}}} \\
=&\varrho(\mathrm{d}(a_1),a_2)u+\varrho(a_1,\mathrm{d}(a_2))u+\varrho(a_1,a_2)\mathrm{d}_\mathfrak{M}(u)+\lambda \varrho(a_1,a_2)u.
\end{align*}
Hence, $(\mathfrak{M},\varrho, \mathrm{d}_\mathfrak{M})$ is a representation over  $(\mathfrak{A}, [-,-,-], \mathrm{d})$.

Since $(\widehat{\mathfrak{A}}, [-,-,-]_{\widehat {\mathfrak{A}}}, \widehat {\mathrm{d}})$ is an abelian extension of $(\mathfrak{A}, [-,-,-], \mathrm{d})$ by  $(\mathfrak{M},[-,-,-]_\mathfrak{M}, \mathrm{d}_\mathfrak{M})$,
by    Proposition~\ref{prop:2-cocycle}, $(\upsilon,\mu)$ is a 2-cocycle.
Moreover, Let $s_1,s_2:\mathfrak{A}\rightarrow \widehat{ \mathfrak{M}}$ be two distinct sections providing 2-cocycles $(\upsilon_1,\mu_1)$ and $(\upsilon_2,\mu_2)$ respectively. Define
linear map $\iota: \mathfrak{A}\rightarrow \mathfrak{M}$ by $\iota(a_1)=s_1(a_1)-s_2(a_1)$. Then
\begin{align*}
&\upsilon_1(a_1,a_2,a_3)\\
=&[s_1(a_1), s_1(a_2),s_1(a_3)]_{\widehat{ \mathfrak{A}}_1}-s_1([a_1,a_2,a_3])\\
=&[s_2(a_1)+\iota(a_1), s_2(a_2)+\iota(a_2), s_2(a_3)+\iota(a_3)]_{\widehat{ \mathfrak{A}}_1}-(s_2([a_1,a_2,a_3])+\iota  ([a_1,a_2,a_3]))\\
=&[s_2(a_1), s_2(a_2), s_2(a_3)]_{\widehat{ \mathfrak{A}}_2}+\varrho(a_2,a_3)\iota(a_1)+\varrho(a_3,a_1)\iota(a_2)+\varrho(a_1,a_2)\iota(a_3)\\
&-s_2([a_1,a_2,a_3])-\iota([a_1,a_2,a_3])\\
=&[s_2(a_1),s_2(a_2),s_2(a_3)]_{\widehat{ \mathfrak{A}}_2}-s_2([a_1,a_2,a_3])\\
&+\varrho(a_2,a_3)\iota(a_1)+\varrho(a_3,a_1)\iota(a_2)+\varrho(a_1,a_2)\iota(a_3)-\iota([a_1,a_2,a_3])\\
=&\upsilon_2(a_1,a_2,a_3)+\delta\iota(a_1,a_2,a_3)
\end{align*}
and
\begin{align*}
\mu_1(a_1)&=\widehat {d}(s_1(a_1))-s_1(\mathrm{d}(a_1))\\
&=\widehat {\mathrm{d}}(s_2(a_1)+\iota(a_1))-\big(s_2(\mathrm{d}(a_1))+\iota(\mathrm{d}(a_1))\big)\\
&=\big(\widehat {\mathrm{d}}(s_2(a_1))-s_2(\mathrm{d}(a_1))\big)+\widehat {\mathrm{d}}(\iota(a_1))-\iota(\mathrm{d}(a_1))\\
&=\mu_2(a_1)+\mathrm{d}_\mathfrak{M}(\iota(a_1))-\iota(\mathrm{d}(a_1))\\
&=\mu_2(a_1)-\Phi\iota(a_1),
\end{align*}
which implies that $(\upsilon_1,\mu_1)-(\upsilon_2,\mu_2)=(\delta\iota,-\Phi\iota)=\partial(\iota)\in \mathcal{C}_{\mathrm{md3Lie_{\lambda}}}^2(\mathfrak{A},\mathfrak{M})$. So $[(\upsilon_1,\mu_1)]=[(\upsilon_2,\mu_2)]\in\mathcal{H}_{\mathrm{md3Lie_{\lambda}}}^2(\mathfrak{A},\mathfrak{M})$.
\end{proof}

 \begin{defn}
Let $(\widehat{\mathfrak{A}}_1, [-,-,-]_{\widehat {\mathfrak{A}}_1}, \widehat {\mathrm{d}}_1)$ and $(\widehat{\mathfrak{A}}_2, [-,-,-]_{\widehat {\mathfrak{A}}_2}, \widehat {\mathrm{d}}_2)$ be two abelian extensions of $(\mathfrak{A}, [-,-,-], \mathrm{d})$ by $(\mathfrak{M},[-,-,-]_\mathfrak{M}, \mathrm{d}_\mathfrak{M})$. They are said to be  equivalent if  there is an isomorphism of modified $\lambda$-differential 3-Lie algebras $\eta:(\widehat{\mathfrak{A}}_1, [-,-,-]_{\widehat {\mathfrak{A}}_1}, \widehat {\mathrm{d}}_1)\rightarrow (\widehat{\mathfrak{A}}_2, [-,-,-]_{\widehat {\mathfrak{A}}_2}, \widehat {\mathrm{d}}_2)$
such that the following diagram is  commutative:
$$\begin{CD}
0@>>> {(\mathfrak{M},\mathrm{d}_\mathfrak{V})} @>i_1 >> (\widehat{\mathfrak{A}}_1,   \widehat {\mathrm{d}}_1) @>p_1 >> (\mathfrak{A},\mathrm{d}) @>>>0\\
@. @| @V \eta VV @| @.\\
0@>>> {(\mathfrak{M},\mathrm{d}_\mathfrak{V})} @>i_2 >> (\widehat{\mathfrak{A}}_2,   \widehat {\mathrm{d}}_2) @>p_2 >> (\mathfrak{A},\mathrm{d}) @>>>0.
\end{CD}$$
\end{defn}

Next we are ready to classify abelian extensions of a modified $\lambda$-differential 3-Lie algebra.

\begin{theorem}
There is a one-to-one correspondence between equivalence classes of abelian extensions of a modified $\lambda$-differential 3-Lie algebra $(\mathfrak{A}, [-,-,-], \mathrm{d})$ by $(\mathfrak{M},[-,-,-]_\mathfrak{M}, \mathrm{d}_\mathfrak{M})$ and the second cohomology group $\mathcal{H}_{\mathrm{md3Lie_{\lambda}}}^2(\mathfrak{A},\mathfrak{M})$ of $(\mathfrak{A}, [-,-,-], \mathrm{d})$ with coefficients in the representation $(\mathfrak{M},\varrho, \mathrm{d}_\mathfrak{M})$.
\end{theorem}
\begin{proof}
  Assume that $(\widehat{\mathfrak{A}}_1, [-,-,-]_{\hat {\mathfrak{A}}_1}, \widehat {\mathrm{d}}_1)$ and $(\widehat{\mathfrak{A}}_2, [-,-,-]_{\hat {\mathfrak{A}}_2}, \widehat {\mathrm{d}}_2)$  are two equivalent abelian extensions of $(\mathfrak{A},$
$ [-,-,-], \mathrm{d})$ by $(\mathfrak{M},[-,-,-]_\mathfrak{M}, \mathrm{d}_\mathfrak{M})$ with the associated isomorphism $\eta: (\widehat{\mathfrak{A}}_1, [-,-,-]_{\hat {\mathfrak{A}}_1}, \widehat {\mathrm{d}}_1) \rightarrow (\widehat{\mathfrak{A}}_2, [-,-,-]_{\hat {\mathfrak{A}}_2}, \widehat {\mathrm{d}}_2)$. Let $s_1$ be a section of $(\widehat{\mathfrak{A}}_1, [-,-,-]_{\hat {\mathfrak{A}}_1}, \widehat {\mathrm{d}}_1)$. As $p_2\circ\eta=p_1$, we have
$$p_2\circ(\eta\circ s_1)=p_1\circ s_1= \mathrm{id}_{\mathfrak{A}}.$$
That is, $\eta\circ s_1$ is a section of $(\widehat{\mathfrak{A}}_2, [-,-,-]_{\hat {\mathfrak{A}}_2}, \widehat {\mathrm{d}}_2)$. Denote $s_2:=\eta\circ s_1$. Since $\eta$ is a isomorphism of  modified $\lambda$-differential 3-Lie algebras such that $\eta|_\mathfrak{M}=\mathrm{id}_\mathfrak{M}$, we get
\begin{align*}
\upsilon_2(a_1,a_2,a_3)&=[s_2(a_1), s_2(a_2), s_{2}(a_3)]_{\widehat{\mathfrak{A}}_2}-s_2([a_1,a_2,a_3])\\
&=[\eta(s_1(a_1)), \eta(s_1(a_2)),\eta(s_1(a_3))]_{\widehat{\mathfrak{A}}_2}-\eta(s_1([a_1,a_2,a_3]))\\
&=\eta\big([s_1(a_1), s_1(a_2), s_{1}(a_3)]_{\widehat{\mathfrak{A}}_1}-s_1([a_1,a_2,a_3])\big)\\
&=\eta(\upsilon_1(a_1,a_2,a_3))\\
&=\upsilon_1(a_1,a_2,a_3)
\end{align*}
and
\begin{align*}
\mu_2(a_1)&=\widehat{\mathrm{d}}_2(s_2(a_1))-s_2(\mathrm{d}(a_1))=\widehat{\mathrm{d}}_2\big(\eta(s_1(a_1))\big)-\eta\big(s_1(\mathrm{d}(a_1))\big)\\
&=\eta\big(\widehat{\mathrm{d}}_1(s_1(a_1))-s_1(\mathrm{d}(a_1))\big)\\
&=\eta(\mu_1(a_1))\\
&=\mu_1(a_1).
\end{align*}
Hence, all equivalent abelian extensions give rise to the same element in $\mathcal{H}_{\mathrm{md3Lie_{\lambda}}}^2(\mathfrak{A},\mathfrak{M})$.

Conversely, suppose that $[(f_1,g_1)]=[(f_2,g_2)]\in\mathcal{H}_{\mathrm{md3Lie_{\lambda}}}^2(\mathfrak{A},\mathfrak{M})$, we can construct two abelian extensions
$0\rightarrow  (\mathfrak{M},[-,-,-]_\mathfrak{M}, \mathrm{d}_\mathfrak{M})\hookrightarrow (\mathfrak{A}\oplus \mathfrak{M},[-,-,-]_{\rho f_1},\mathrm{d}_{g_1})\stackrel{p_1}{\longrightarrow} (\mathfrak{A}, [-,-,-], \mathrm{d})\rightarrow 0$ and  $0\rightarrow  (\mathfrak{M},[-,-,-]_\mathfrak{M}, \mathrm{d}_\mathfrak{M})\hookrightarrow (\mathfrak{A}\oplus \mathfrak{M},[-,-,-]_{\rho f_2},\mathrm{d}_{g_2})\stackrel{p_2}{\longrightarrow} (\mathfrak{A}, [-,-,-], \mathrm{d})\rightarrow 0$ via Eqs.~\eqref{4.1} and \eqref{4.2}. Then there exists  a linear map $\iota: \mathfrak{A}\rightarrow  \mathfrak{M}$ such that
 $$(f_2,g_2)=(f_1,g_1)+\partial(\iota).$$
 Define linear map $\eta_\iota: \mathfrak{A}\oplus \mathfrak{M}\rightarrow  \mathfrak{A}\oplus \mathfrak{M}$ by
$\eta_\iota(a_1,u_1):=a_1+\iota(a_1)+u_1, ~a_1\in \mathfrak{A}, u_1\in \mathfrak{M}.$
Then, $\eta_\iota$ is an isomorphism of these two abelian extensions $(\mathfrak{A}\oplus \mathfrak{M},[-,-,-]_{\rho f_1},\mathrm{d}_{g_1})$ and $(\mathfrak{A}\oplus \mathfrak{M},[-,-,-]_{\rho f_2},\mathrm{d}_{g_2})$.
\end{proof}

\begin{remark}
In particular, any vector space $\mathfrak{M}$ with linear transformation $\mathrm{d}_\mathfrak{M}$ can serve as a trivial
representation of $(\mathfrak{A}, [-,-,-], \mathrm{d})$. In this situation, central extensions of $(\mathfrak{A}, [-,-,-], \mathrm{d})$ by $(\mathfrak{M},[-,-,-]_\mathfrak{M}, \mathrm{d}_\mathfrak{M})$ are classified by
the second cohomology group $\mathcal{H}_{\mathrm{md3Lie_{\lambda}}}^2(\mathfrak{A},\mathfrak{M})$ of $(\mathfrak{A}, [-,-,-], \mathrm{d})$ with the coefficient in the trivial representation
 $(\mathfrak{M}, \rho=0,  \mathrm{d}_\mathfrak{M}).$
\end{remark}

\section{$T^*$-extensions of modified $\lambda$-differential 3-Lie algebras} \label{sec:Text}
\def\theequation{\arabic{section}.\arabic{equation}}
\setcounter{equation} {0}
The $T^*$-extensions of  3-Lie algebra was studied in \cite{LiuW}. In this section, we consider $T^*$-extensions of modified $\lambda$-differential 3-Lie algebras  by the second cohomology groups with the coefficient  in a dual adjoint representation.

Let $(\mathfrak{A}, [-, -, -],\mathrm{d})$ be a   modified $\lambda$-differential 3-Lie algebra and $\mathfrak{A}^*$ be the dual space  of  $\mathfrak{A}$. By   Example~\ref{exam:dual rep},
$(\mathfrak{A}^*;    ad^* , -\mathrm{d}^*)$ is a   representation of  $(\mathfrak{A}, [-, -, -],\mathrm{d})$.
Suppose that $(f,g)\in\mathcal{C}_{\mathrm{md3Lie_\lambda}}^2(\mathfrak{A}, \mathfrak{A}^*)$.
Define a trilinear map  $[-,-,-]_{f}: \wedge^3 (\mathfrak{A}\oplus\mathfrak{A}^*)\rightarrow \mathfrak{A}\oplus\mathfrak{A}^*$
and  a linear map   $ \mathrm{d}_{g}: \mathfrak{A}\oplus\mathfrak{A}^*\rightarrow \mathfrak{A}\oplus\mathfrak{A}^*$ respectively by
\begin{align}
&[a_1+\alpha_1, a_2+\alpha_2,a_3+\alpha_3]_{ f} \nonumber\\
&=[a_1, a_2, a_3]+ad^*(a_2,a_3)\alpha_1+ad^*(a_3,a_1)\alpha_2+ad^*(a_1,a_2)\alpha_3+f(a_1, a_2, a_3), \label{5.1} \\
&\mathrm{d}_g(a_1+\alpha_1)=\mathrm{d}(a_1)-\mathrm{d}^*(\alpha_1)+g(a_1),\,\forall a_1, a_2, a_3\in \mathfrak{A},\ \alpha_1, \alpha_2, \alpha_3\in \mathfrak{A}^*.\label{5.2}
\end{align}

Similar to Proposition \ref{prop:2-cocycle}, we have the following  result.

\begin{prop} \label{prop:2-cocyclead}
 With the above notations,  $(\mathfrak{A}\oplus \mathfrak{A}^*,[-,-,-]_{f},\mathrm{d}_{g})$ is a modified $\lambda$-differential 3-Lie algebra  if and only if
$(f,g)$ is a 2-cocycle   in the cohomology of the modified $\lambda$-differential 3-Lie algebra $(\mathfrak{A}, [-, -, -],\mathrm{d})$ with the coefficient  in the representation $(\mathfrak{A}^*;  ad^*, -\mathrm{d}^*)$.
\end{prop}

\begin{defn}
The modified $\lambda$-differential 3-Lie algebra $(\mathfrak{A}\oplus \mathfrak{A}^*,[-,-,-]_{f},\mathrm{d}_{g})$ is called the $T^*$-extension of the
modified $\lambda$-differential 3-Lie algebra $(\mathfrak{A}, [-, -, -],\mathrm{d})$. Denote the $T^*$-extension by  $T^*_{(f,g)}(\mathfrak{A})=(T^*(\mathfrak{A})=\mathfrak{A}\oplus \mathfrak{A}^*,[-,-,-]_{f},\mathrm{d}_{g})$.
\end{defn}

\begin{defn}
Let $(\mathfrak{A}, [-, -, -],\mathrm{d})$ be a   modified $\lambda$-differential 3-Lie algebra. $(\mathfrak{A}, [-, -, -],\mathrm{d})$ is said to be
metrised if it has a non-degenerate symmetric bilinear form $\varpi_\mathfrak{A}$ which satisfying
\begin{align}
&\varpi_\mathfrak{A}([a_1,a_2,a_3],a_4)+\varpi_\mathfrak{A}(a_3,[a_1,a_2,a_4])=0, \label{5.3} \\
&\varpi_\mathfrak{A}(\mathrm{d}(a_1),a_2)+\varpi_\mathfrak{A}(a_1,\mathrm{d}(a_2))=0, \ \  \forall a_1,a_2,a_3,a_4\in\mathfrak{A}.\label{5.4}
\end{align}
We may also say that $(\mathfrak{A}, [-, -, -],\mathrm{d},\varpi_\mathfrak{A})$ is a metrised modified $\lambda$-differential 3-Lie algebra.
\end{defn}
Define a bilinear map $\varpi:\wedge^2T^*(\mathfrak{A})\rightarrow \mathfrak{A}$ by
\begin{align}
&\varpi(a_1+\alpha_1,a_2+\alpha_2)=\alpha_1(a_2)+\alpha_2(a_1), \  \ \forall a_1, a_2\in \mathfrak{A}, \alpha_1, \alpha_2\in \mathfrak{A}^* \label{5.5}
\end{align}

\begin{prop} \label{prop:metrised}
 With the above notations,  $(T^*_{(f,g)}(\mathfrak{A}), \varpi)$ is a metrised modified $\lambda$-differential 3-Lie algebra  if and only if
$$f(a_1,a_2,a_3)(a_4)+f(a_1,a_2,a_4)(a_3)=0, \ \ g(a_1)(a_2)+g(a_2)(a_1)=0, \ \ \forall a_1,a_2,a_3,a_4\in \mathfrak{A}.$$
\end{prop}
\begin{proof}
For any $a_1, a_2, a_3,a_4\in \mathfrak{A},\ \alpha_1, \alpha_2, \alpha_3,\alpha_4\in \mathfrak{A}^*$, using Eqs.~\eqref{2.6},  \eqref{5.1}-\eqref{5.5} we have
\begin{align*}
&\varpi([a_1+\alpha_1,a_2+\alpha_2,a_3+\alpha_3]_f,a_4+\alpha_4)+\varpi(a_3+\alpha_3,[a_1+\alpha_1,a_2+\alpha_2,a_4+\alpha_4]_f)\\
=&\varpi([a_1, a_2, a_3]+ad^*(a_2,a_3)\alpha_1+ad^*(a_3,a_1)\alpha_2+ad^*(a_1,a_2)\alpha_3+f(a_1, a_2, a_3),a_4+\alpha_4)\\
&+\varpi(a_3+\alpha_3,[a_1, a_2, a_4]+ad^*(a_2,a_4)\alpha_1+ad^*(a_4,a_1)\alpha_2+ad^*(a_1,a_2)\alpha_4+f(a_1, a_2, a_4))\\
=&\alpha_4([a_1, a_2, a_3])+ad^*(a_2,a_3)\alpha_1(a_4)+ad^*(a_3,a_1)\alpha_2(a_4)+ad^*(a_1,a_2)\alpha_3(a_4)+f(a_1, a_2, a_3)(a_4)\\
&+\alpha_3([a_1, a_2, a_4])+ad^*(a_2,a_4)\alpha_1(a_3)+ad^*(a_4,a_1)\alpha_2(a_3)+ad^*(a_1,a_2)\alpha_4(a_3)+f(a_1, a_2, a_4)(a_3)\\
=&\alpha_4([a_1, a_2, a_3])-\alpha_1 ([a_2,a_3,a_4]-\alpha_2([a_3,a_1,a_4])-\alpha_3([a_1,a_2,a_4])+f(a_1, a_2, a_3)(a_4)\\
&+\alpha_3([a_1, a_2, a_4])-\alpha_1 ([a_2,a_4,a_3])-\alpha_2 ([a_4,a_1,a_3])-\alpha_4(a_1,a_2,a_3)+f(a_1, a_2, a_4)(a_3)\\
=&f(a_1, a_2, a_3)(a_4)+f(a_1, a_2, a_4)(a_3)\\
=&0,\\
&\varpi(\mathrm{d}_g(a_1+\alpha_1),a_2+\alpha_2)+\varpi(a_1+\alpha_1,\mathrm{d}_g(a_2+\alpha_2))\\
=&\varpi(\mathrm{d}(a_1)-\mathrm{d}^*(\alpha_1)+g(a_1),a_2+\alpha_2)+\varpi(a_1+\alpha_1,\mathrm{d}(a_2)-\mathrm{d}^*(\alpha_2)+g(a_2))\\
=&-\mathrm{d}^*(\alpha_1)(a_2)+g(a_1)(a_2)+\alpha_2(\mathrm{d}(a_1))+\alpha_1(\mathrm{d}(a_2))-\mathrm{d}^*(\alpha_2)(a_1)+g(a_2)(a_1)\\
=&-\alpha_1(\mathrm{d}(a_2))+g(a_1)(a_2)+\alpha_2(\mathrm{d}(a_1))+\alpha_1(\mathrm{d}(a_2))-\alpha_2(\mathrm{d}(a_1))+g(a_2)(a_1)\\
=&g(a_1)(a_2)+g(a_2)(a_1)\\
=&0.
\end{align*}
Thus,  we get the result.
\end{proof}

\noindent
{{\bf Acknowledgments.}  The paper is  supported  by  the Science and Technology Program of Guizhou Province (Grant Nos. QKHZC[2023]372),
the Scientific Research Foundation for Science \& Technology Innovation Talent Team of the Intelligent Computing and Monitoring of Guizhou Province (Grant No. QJJ[2023]063),   the National Natural Science Foundation of China (Grant No. 12161013).

\end{document}